\documentclass[journal]{IEEEtran}
\usepackage{comment}
\usepackage{soul}
\usepackage{color}
\usepackage[usenames,dvipsnames]{xcolor}

\usepackage{cite}
\usepackage{url}

\ifCLASSINFOpdf
\else
\fi
\usepackage{tikz}
\usetikzlibrary{arrows, positioning, shapes.geometric, fit, calc,backgrounds,shapes.arrows}
\usetikzlibrary{arrows.meta}

\usepackage{amsmath}
\usepackage[short]{optidef}
\usepackage[ruled,vlined]{algorithm2e}

\ifCLASSOPTIONcompsoc
  \usepackage[caption=false,font=normalsize,labelfont=sf,textfont=sf]{subfig}
\else
  \usepackage[caption=false,font=footnotesize]{subfig}
\fi

\renewcommand{\baselinestretch}{0.96}

\usepackage{bm}
\usepackage{hyperref}
\usepackage[capitalise]{cleveref}
\Crefname{equation}{}{}

\begin{document}

\title{Design of a Continuous Local Flexibility Market with Network Constraints \\
\thanks{This work is supported by the H2020 European Project FLEXGRID, Grant Agreement No. 863876.}
}

\author{\IEEEauthorblockN{Eléa~Prat, Lars~Herre, Jalal Kazempour,  Spyros~Chatzivasileiadis} \\
\IEEEauthorblockA{Department of Electrical Engineering, Technical University of Denmark, Kongens Lyngby, Denmark} \\
\{emapr, lfihe, seykaz, spchatz\}@elektro.dtu.dk
\vspace{-0.7cm}
}

\maketitle


\begin{abstract}
To the best of our knowledge, this paper proposes for the first time a design of a continuous local flexibility market that explicitly considers network constraints. Continuous markets are expected to be the most appropriate design option during the early stages of local flexibility markets, where insufficient liquidity can hinder market development. At the same time, increasingly loaded distribution systems require to explicitly consider network constraints in local flexibility market clearing in order to help resolve rather than aggravate local network problems, such as line congestion and voltage issues. This paper defines the essential design considerations, introduces the local flexibility market clearing algorithm, and -- aiming to establish a starting point for future research -- discusses design options and research challenges that emerge during this procedure which require further investigation. 
\end{abstract}


%
\IEEEpeerreviewmaketitle

\vspace{-0.1cm}
\section{Introduction}
The proliferation of energy resources with variable and uncertain power profiles at the low-voltage distribution level (e.g., renewables, electric vehicles) call for drastically higher levels of flexibility. 
The concept of \textit{local flexibility markets} at the distribution level has emerged recently \cite{Jin2020}, and is investigated in several EU projects\footnote{See \url{http://www.interrface.eu/}, \url{http://smartnet-project.eu/}, \mbox{\url{https://flexgrid-project.eu/}}, and \url{http://www.eu-ecogrid.net/}.}, such as INTERRFACE, SmartNet, and FLEXGRID, as well as national projects such as EcoGrid 2.0 in Denmark.
 Local flexibility markets are expected to increase the reliability in the power supply, and, at the same time, help avoid local problems such as line congestion and voltage issues in the distribution network. Considering the increasingly loaded distribution systems, incorporating the network constraints in such market clearing algorithms is necessary, so that the procured flexibility helps resolve and not further aggravate existing network problems.
 In an envisioned flexibility market, the flexibility product may be traded in the form of energy, balancing capacity or reserve capacity.
This paper designs a local market to trade \textit{flexibility reserve capacity} in a forward stage while ensuring operational feasibility of the real-time activation. 

Usual modeling approaches either ignore the network \cite{Correa2020, Lampropoulos2019, Ziras2021,Olivella-Rosell2018}, or if they include the local flexibility scheduling in distribution power flow calculations, they assume that the distribution system operator (DSO) and the flexibility market operator (FMO) form one entity \cite{Morstyn2019,Alanazi2020,Pastor2018, Torbaghan2020}. In that case, decisions about flexibility procurement and activation are made considering requirements of DSO to solve voltage or congestion issues. %
However, the legal framework may (and, in the EU, currently does) not allow the DSO to act as the market operator. 
On the contrary, Ref.~\cite{Heinrich2020} explicitly models the FMO as a separate entity which requires access to distribution network data in order to ensure operational feasibility, while also explicitly considers network constraints. However, their approach relies on an auction-based market clearing algorithm. 

To the best of our knowledge, and as summarized in \cref{tab:Literature}, so far no work exists that has proposed a design for a \emph{continuous} local flexibility market which also includes the network constraints. 

In contrast to auctions that close once or on multiple sequential gate closures, a continuous market clears as soon as a pair of bids matches. On the other hand, auctions allow for different pricing mechanisms, e.g., uniform,  pay-as-bid or Vickrey–Clarke–Groves, while the continuous market requires a pay-as-bid pricing to be used \cite{Schittekatte2020} in practice.
Pilot projects, such as Piclo Flex, use an auction-based market design similar to those in the wholesale energy markets, whereas Enera, GOPACS, and NODES\footnote{See \url{https://piclo.energy/}, \url{https://projekt-enera.de/}, \url{https://gopacs.eu/}, and \url{https://nodesmarket.com/}.} implement continuous trading with pay-as-bid pricing (without, however, explicitly considering the network constraints in the market clearing).
Ref.~\cite{Lauterbach1997} suggests that continuous trading promotes price efficiency and better suits the trading mechanism preferences of the investors, while Ref.~\cite{Schittekatte2020} suggests that continuous markets might be more suitable for markets with lower liquidity. In the early stages of local flexibility markets, where insufficient liquidity may hinder market development, we expect that continuous markets would be the appropriate design option, as also suggested by the several pilot projects. 

\begin{table}[t]
\centering
\caption{Modeling approaches in the related literature}
\begin{tabular}{cccc}
\hline
Ref.               & Trading type   & FMO      & Network check\\ \hline
\cite{Morstyn2019}      & Negotiation & One entity with DSO   & Yes (AC-OPF) \\
\cite{Alanazi2020,Pastor2018}      & Auction   & One entity with DSO     & Yes (AC-OPF$^{1}$)\\
\cite{Torbaghan2020}    & Auction   & One entity with DSO       & Yes (SOC$^{2}$)\\
\cite{Olivella-Rosell2018,Ziras2021}& Auction & One entity with DSO     & No \\
\cite{Lampropoulos2019,Correa2020} & Auction   & Separate entity & No \\
\cite{Heinrich2020}     & Auction   & Separate entity & Yes$^{3}$\\
\hline
This paper              & Continuous& Separate entity & Yes (DC) \\
\hline
\multicolumn{4}{l}{\tiny{$^{1}$ In \cite{Alanazi2020}, the nonlinear AC optimal power flow (OPF) equations are linearized.}}\\
\multicolumn{4}{l}{\tiny{$^{2}$ \cite{Torbaghan2020} uses a second order cone (SOC) relaxation of optimal flexibility dispatch including line and voltage constraints.}}\\
\multicolumn{4}{l}{\tiny{$^{3}$ It is unclear which power flow algorithm is used in \cite{Heinrich2020}.}}
\end{tabular}
\label{tab:Literature}
\vspace{-0.5cm}
\end{table}

To the best of our knowledge, this is the first paper that proposes the design of a continuous local flexibility market that explicitly considers network constraints. The contributions of this paper are the following:
\begin{itemize}
    \item We introduce a continuous market clearing mechanism for local flexibility markets which considers network constraints. The focus of this paper is on active power markets and considers active power flows and line limits.
    \item We define the essential design considerations and discuss challenges and design options that arise during the design of such a local flexibility market.
    We also suggest directions that require further research.
\end{itemize}

The rest of this paper is structured as follows.
In \cref{sec:MarketArch} we detail the network-aware continuous market clearing architecture of the local flexibility market operator. \cref{sec:Case} applies the designed local flexibility market on a distribution system. In \cref{sec:Disc}, we discuss the implications of the proposed design, as highlighted in the case study, and suggest future directions for research. \cref{sec:Conc} concludes.

\vspace{-0.1cm}
\section{Market Architecture} \label{sec:MarketArch}

\subsection{General Characteristics} \label{sec:Charac}

The market we envision is a distribution-level local flexibility market that clears continuously with the pay-as-bid pricing rule. Market actors include the DSO as well as balance responsible parties (BRPs).
Actors submit a bid as \emph{FlexRequest} or \emph{FlexOffer} for active power reserve capacity (availability) in either upward or downward direction. Here, upward indicates an increase of production or a decrease in consumption, while downward indicates a decrease in production or an increase in consumption. The bid is composed of its type, i.e., \emph{FlexOffer} or \emph{FlexRequest} and \emph{up} or \emph{down}, price, volume, and location (network bus). 
Incoming, non-matching bids are placed in the \textit{order book} until they are cleared with a matching bid. \emph{FlexRequest} and \emph{FlexOffer} are cleared with a set of rules described in the following.

The first-come first-served principle is used to match a \textit{FlexOffer} with a \textit{FlexRequest}, where the price is set by the bid that came in first. 
Other options for pricing mechanisms are discussed in \cref{sec:Disc}.
If the volume and price allow \textit{FlexOffer} and \textit{FlexRequest} to match, a network check is performed by the FMO before the bids can be cleared. The location of the bid does not need to match, i.e., \emph{FlexOffer} and \emph{FlexRequest} can be located at different buses. 
The network check is based on a baseline energy dispatch that is established by either previous markets (e.g., day-ahead energy market) or by an estimation of load and generation at each bus (based on, e.g., usually available data of similar days and hours and load forecasting).
Since DSOs are, e.g., in the EU, currently not allowed to act as FMO, the network-aware market clearing requires data exchange between DSO and FMO. Specifically, the DSO needs to share the network data with the FMO, similar to the data exchange on transmission level.

Owing to this network check requirement, it is especially important to allow partial matching of the bids. In this way, we can make sure that two bids can match up to the point where their activation could result in a congestion.

The market architecture introduced here is general enough to be integrated in any framework. In particular, time structures (gate closure time, market resolution) and interaction with the wholesale markets are not discussed in this work.

\vspace{-0.1cm}
\subsection{Flexibility Requests (FlexRequest)}
 In a first step, the DSO would be the main buyer of flexibility. In this perspective, the DSO would formulate a \emph{FlexRequest} as an abstract representation of a contingency which is submitted to the FMO in the form of a bid.
 In order to formulate a general approach that also holds for future players and needs, we do not limit our framework to the DSO alone.
A fundamental question we have raised is whether the bid should be required to include a location or not. Indeed, the flexibility buyer could be estimating an aggregated need for upward or downward flexibility, without knowing or desiring to share how this is going to be split among the different buses.
However, a network check algorithm (see \cref{sec:NetworkCheck} for more details) requires to associate a location with every \textit{FlexRequest} (and \textit{FlexOffer}) in order to be able to assess their feasibility if they match. If the \textit{FlexRequest} does not determine location, the algorithm shall check for all possible locations of the \textit{FlexRequest}, and ensure that the matching bids remain feasible for all. Besides increasing complexity exponentially (combinatorial problem), this, most importantly, decreases the chances for bids to match, as the probability of finding one of the potential location pairs infeasible increases substantially. In our experiments, we found out that a market-clearing algorithm that does not require to determine location for \textit{FlexRequests} would successfully match bids in a significantly unconstrained network. Such a network, however, resembles a copperplate. In that case, a network check would be redundant, and the flexibility market clearing could use existing standard schemes. Our goal in this paper is to design a flexibility market that can be used in constrained distribution networks and allow bids to match only if they do not lead to any congestion or violate network constraints. Therefore, in the design proposed in this paper, \textit{FlexRequests} shall include the location. In future work, we plan to relax this requirement and explore whether the design of scenarios about the most probable locations for the submitted \textit{FlexRequests} can lead to an efficient network-aware market which will be feasible with high probability. 

In addition to the inclusion or not of the location, a feature that we have included on the design of the proposed market is that the \textit{FlexRequests} can specify whether they are \textit{conditional} or \textit{unconditional}. Based on prior work \cite{corey}, market actors seem to be in a position to estimate whether their \textit{FlexRequest} will be activated in the real-time operation with high probability (certainty) or not. A request tagged as unconditional is expected to be activated with certainty, unlike a request tagged as conditional. With this feature, one can consider that the market can be used both to clear energy and capacity.

\vspace{-0.1cm}
\subsection{Network Check: Insights from a 3-Bus System} \label{sec:NetworkCheck}
The state-of-the-art in reserve clearing on transmission level is a security-constrained DC-OPF that ensures that for each contingency individually, there is a feasible combination of reserve activation. Here, we expand on this approach by also considering combinations of contingencies which are represented by \emph{FlexRequests}.

\subsubsection{Network Model}
Regarding the network check, the first decision to make is which power flow algorithm to use. In this paper, we have used the DC power flow algorithm as the first step towards the inclusion of network constraints in a continuous market clearing algorithm. As we elaborate in \cref{sec:Disc}, two main reasons for this choice is that the DC power flow is simpler, and thus more transparent for the market players, and faster, with computing time being a critical element for continuous markets (please see Section~IVa for more details). Future work will include the extension of this algorithm to LinDistFlow \cite{BaranWu1989_LinDistFlow} and AC power flow.

\subsubsection{Check Procedure}
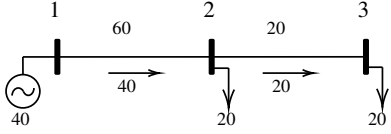
\begin{figure}[t]
  \centering 
  \resizebox{0.60\columnwidth}{!}{\tikzset{every picture/.style={line width=0.75pt}} 

\begin{tikzpicture}[x=0.75pt,y=0.75pt,yscale=-1,xscale=1]

\draw    (40.14,140.07) -- (241.25,140) ;
\draw  [fill={rgb, 255:red, 0; green, 0; blue, 0 }  ,fill opacity=1 ] (59,130.5) .. controls (59,130.22) and (59.22,130) .. (59.5,130) -- (61,130) .. controls (61.28,130) and (61.5,130.22) .. (61.5,130.5) -- (61.5,149.5) .. controls (61.5,149.78) and (61.28,150) .. (61,150) -- (59.5,150) .. controls (59.22,150) and (59,149.78) .. (59,149.5) -- cycle ;
\draw  [fill={rgb, 255:red, 0; green, 0; blue, 0 }  ,fill opacity=1 ] (149,130.5) .. controls (149,130.22) and (149.22,130) .. (149.5,130) -- (151,130) .. controls (151.28,130) and (151.5,130.22) .. (151.5,130.5) -- (151.5,149.5) .. controls (151.5,149.78) and (151.28,150) .. (151,150) -- (149.5,150) .. controls (149.22,150) and (149,149.78) .. (149,149.5) -- cycle ;
\draw  [fill={rgb, 255:red, 0; green, 0; blue, 0 }  ,fill opacity=1 ] (239,130.5) .. controls (239,130.22) and (239.22,130) .. (239.5,130) -- (241,130) .. controls (241.28,130) and (241.5,130.22) .. (241.5,130.5) -- (241.5,149.5) .. controls (241.5,149.78) and (241.28,150) .. (241,150) -- (239.5,150) .. controls (239.22,150) and (239,149.78) .. (239,149.5) -- cycle ;
\draw    (160.14,146.55) -- (160.14,168.07) ;
\draw [shift={(160.14,170.07)}, rotate = 270] [color={rgb, 255:red, 0; green, 0; blue, 0 }  ][line width=0.75]    (10.93,-3.29) .. controls (6.95,-1.4) and (3.31,-0.3) .. (0,0) .. controls (3.31,0.3) and (6.95,1.4) .. (10.93,3.29)   ;
\draw    (90,149.5) -- (117.86,149.51) ;
\draw [shift={(119.86,149.51)}, rotate = 180.03] [color={rgb, 255:red, 0; green, 0; blue, 0 }  ][line width=0.75]    (6.56,-1.97) .. controls (4.17,-0.84) and (1.99,-0.18) .. (0,0) .. controls (1.99,0.18) and (4.17,0.84) .. (6.56,1.97)   ;
\draw    (180,149.5) -- (207.86,149.51) ;
\draw [shift={(209.86,149.51)}, rotate = 180.03] [color={rgb, 255:red, 0; green, 0; blue, 0 }  ][line width=0.75]    (6.56,-1.97) .. controls (4.17,-0.84) and (1.99,-0.18) .. (0,0) .. controls (1.99,0.18) and (4.17,0.84) .. (6.56,1.97)   ;
\draw   (30.32,160.11) .. controls (30.32,154.74) and (34.73,150.39) .. (40.18,150.39) .. controls (45.62,150.39) and (50.04,154.74) .. (50.04,160.11) .. controls (50.04,165.47) and (45.62,169.82) .. (40.18,169.82) .. controls (34.73,169.82) and (30.32,165.47) .. (30.32,160.11) -- cycle ;
\draw  [draw opacity=0] (33.84,162.62) .. controls (33.87,160.17) and (35.48,158.2) .. (37.45,158.2) .. controls (38.65,158.2) and (39.72,158.93) .. (40.37,160.06) -- (37.45,162.68) -- cycle ; \draw   (33.84,162.62) .. controls (33.87,160.17) and (35.48,158.2) .. (37.45,158.2) .. controls (38.65,158.2) and (39.72,158.93) .. (40.37,160.06) ;
\draw  [draw opacity=0] (46.6,158.21) .. controls (46.58,160.59) and (45.12,162.5) .. (43.33,162.5) .. controls (42.04,162.5) and (40.94,161.53) .. (40.4,160.1) -- (43.33,158.15) -- cycle ; \draw   (46.6,158.21) .. controls (46.58,160.59) and (45.12,162.5) .. (43.33,162.5) .. controls (42.04,162.5) and (40.94,161.53) .. (40.4,160.1) ;

\draw    (150.14,146.64) -- (160.14,146.64) ;
\draw    (249.86,146.07) -- (249.86,167.6) ;
\draw [shift={(249.86,169.6)}, rotate = 270] [color={rgb, 255:red, 0; green, 0; blue, 0 }  ][line width=0.75]    (10.93,-3.29) .. controls (6.95,-1.4) and (3.31,-0.3) .. (0,0) .. controls (3.31,0.3) and (6.95,1.4) .. (10.93,3.29)   ;
\draw    (239.86,146.07) -- (249.86,146.07) ;
\draw    (40.18,150.39) -- (40.14,140.07) ;

\draw (54,106) node [anchor=north west][inner sep=0.75pt] [align=left] {1};
\draw (144,106) node [anchor=north west][inner sep=0.75pt] [align=left] {2};
\draw (234,106) node [anchor=north west][inner sep=0.75pt] [align=left] {3};
\draw (91,118) node [anchor=north west][inner sep=0.75pt]  [font=\footnotesize] [align=left] {60};
\draw (181,118) node [anchor=north west][inner sep=0.75pt]  [font=\footnotesize] [align=left] {20};
\draw (32,171) node [anchor=north west][inner sep=0.75pt]  [font=\footnotesize] [align=left] {40};
\draw (152,171) node [anchor=north west][inner sep=0.75pt]  [font=\footnotesize] [align=left] {20};
\draw (94,152) node [anchor=north west][inner sep=0.75pt]  [font=\footnotesize] [align=left] {40};
\draw (184,152) node [anchor=north west][inner sep=0.75pt]  [font=\footnotesize] [align=left] {20};
\draw (240,171) node [anchor=north west][inner sep=0.75pt]  [font=\footnotesize] [align=left] {20};

\end{tikzpicture}} 
  \caption{Illustrative 3-bus network, all values in kW. Values on top show the line limits. Values below the arrows show the baseline dispatch.}
  \label{fig:network3bus}
\vspace{-0.1cm}
\end{figure}

When designing the network check algorithm, one has to keep in mind that this is a market for flexibility reserves. There is no guarantee that the procured reserves will be activated, but we need to make sure that they can be activated without causing any congestion. Here, we discuss how to achieve feasible solutions at both the market clearing stage and during real-time activation.
The example of a simple 3-bus system with DC power flow will be used as an illustration. The initial state available to the FMO is shown in Fig. \ref{fig:network3bus}.

The first point to consider is the difference between the quantity procured and the quantity activated. We need to make sure that any activation in the range between zero and the procured quantity would not create any congestion. In particular, one has to keep in mind that it might not be enough to only verify the feasibility of full activation. In our case study in \cref{sec:Case}, we will prove that for the chosen setup, it is. This proof can also be applied to the example given here, so for the rest of the section, we consider the activation of only the complete procured capacities.

The second point to discuss is how to take into account the previous matches between requests and offers when checking the feasibility of the current match. A number of options are possible, which we discuss in the rest of this section. Note that the following only applies to conditional requests, as unconditional requests are considered to be activated in any case. As a consequence, a match with an unconditional request would directly modify the initial dispatch and be used as the new baseline for subsequent bids.

\paragraph{Individual Effect}
The first option would be to ensure that each new match of a request and an offer does not cause any congestion when it is the only one activated. In that case, in our example, both bids shown in \cref{tab:bid1} would be accepted. However, if they were activated together, there would be congestion in line 1-2 and the activation of the second bid would lead to a network constraint violation. As a result, the procured flexibility would either not be delivered or the system would be at risk. It becomes clear that considering only the individual effect of the bids during the network check is a very limited approach, as with high probability more than one conditional bids will be activated at the same time.

\begin{table}[t]
\centering
\caption{Example of accepted bids for the 3-bus system, considering individual effect}
\begin{tabular}{ccc}
\hline
\multicolumn{1}{l}{Bidding  round} & Request  (kW)    & Accepted offer  (kW) \\ \hline
1                                  & 10 downward for 1 & -10 in 2              \\
2                                  & 20 upward for 2   & +20 in 1              \\ \hline
\end{tabular}
\label{tab:bid1}
\end{table}

\paragraph{Cumulative Effect}
A different option would be to consider the effect of all previously accepted bids. Following this procedure, the bids in Table \ref{tab:bid2} would all be accepted. In this example, considering the cumulative effect, we ensure that the bids of the first round can be activated alone, the bids of the first and second rounds can be activated together and all three bids can be activated at the same time. 

\begin{table}[t]
\centering
\caption{Example of accepted bids for the 3-bus system, considering cumulative effect}
\begin{tabular}{ccc}
\hline
\multicolumn{1}{l}{Bidding  round} & Request  (kW.)    & Accepted offer  (kW) \\ \hline
1                                  & 20 upward for 1   & +20 in 2              \\
2                                  & 30 downward for 3 & -30 in 1              \\
3                                  & 20 downward for 2 & -20 in 3              \\ \hline
\end{tabular}
\label{tab:bid2}
\end{table}

In this case, however, note that the bids from the third round cannot be activated alone, because this would lead to the congestion of line 2-3 and the dispatch would be infeasible. However, if we could activate the first two requests, this would remove the congestion, and the third request could then be served. For this to happen, we would need an actor that has access to all matched \textit{FlexOffers} and has the ability to activate them if necessary. This could be the role of the DSO.

\paragraph{Individual and Cumulative Effects}
Checking both individual and cumulative effects for each bid would considerably reduce the cases where the activation is operationally infeasible, but not remove them completely, as the activation of only a subset of bids would not be explicitly considered.

\paragraph{All Combinations}
The only way to make sure that the activation would not lead to network violations (line limit violations in our case) is to test the activation of all combinations of accepted bids with the new bid under check. The issue with this approach is that it results in higher computing times as the number of accepted bids increases, which is critical for a continuous market-clearing algorithm. To reduce that burden, the checks of the different combinations can be easily performed in parallel in this case.

\paragraph{Scenarios}
Depending on the general context, one can decide that it is not necessary to make sure that all possible activation combinations are feasible. Instead, one could consider a set of most probable scenarios for their activation. This could work well in the case that the DSO has access to other solutions to avoid congestion if, for example, a bid activation occurs that was not captured in the scenarios. 
The DSO could give an instruction on the maximum probability that the activation of bids leads to a congestion.

\paragraph{Unconditional Requests}
The last point to discuss is the effect of unconditional requests matching on previously rejected bids. We can get insights from the 3-bus system. If we assume that only unconditional requests were submitted and that all these requests were submitted before the offers are added, the situation depicted in Table \ref{tab:bid3} could arise.

\begin{table}[t]
\centering
\caption{Example of accepted bids for the 3-bus system, with unconditional requests}
\begin{tabular}{ccc}
\hline
\multicolumn{1}{l}{Bidding round} & Submitted offer (kW)    & Matching request (kW) \\ \hline
1                                 & -20 in 3                 & No match (congestion)  \\
2                                 & +30 in 3                 & 30 upward for 1        \\
2                                 & -20 in 3 (re-evaluation) & 20 downward for 2      \\ \hline
\end{tabular}
\label{tab:bid3}
\end{table}

In this case, there are two requests, one for upward flexibility in bus 1 and one for downward flexibility in bus 2. When the first offer is submitted, it cannot match any of the requests because of congestion; so, the offer is added to the order book. But the match with the second offer relieves this congestion. As a consequence, it is important to make sure that the bids in the order book are re-evaluated once unconditional requests are matched, as they modify the power dispatch.

\vspace{-0.1cm}
\section{Case Study}\label{sec:Case}

In this section, an example of such a continuous flexibility market is described. First, the characteristics chosen for the design of this market are given. It is then applied to the 15-bus system from \cite{Das1995}.  The data used and the code for the matching algorithm are available online \cite{GitHub}.

\vspace{-0.1cm}
\subsection{Market Features and Assumptions}
For the following case studies, a simple market setup is assumed. We apply the general characteristics described in Section \ref{sec:MarketArch}. The offers and requests are for active power only. We assume that the initial dispatch is feasible (no congestion).
We also assume that, as a first version of a market, there are no block offers and it is not possible to combine up and down offers for a given request. The power flows are calculated with the help of the power transfer distribution factors (PTDFs), assuming a DC, i.e., linearized, power flow.

When performing the network check, all combinations of the previous matches with the bid-match under check are considered, to ensure that no congestion could result from their activation in real time (individual, all, or a subset of them). Partial match is allowed, following what was discussed in Section \ref{sec:Charac}. When unconditional requests are accepted, all the offers in the order book are evaluated again.

\vspace{-0.1cm}
\subsection{Network Check with PTDFs}

PTDFs are linear sensitivities linking power injections with line flows (for more details, see \cite{notes}). In particular, the power flow in the line between bus $i$ and $j$, $P_{ij}$, is linked to the power injected at bus $m$, $P_m$, by the PTDF factor of line $ij$ for an injection of power at the slack bus $k$ and retrieval of the same quantity in bus $m$, $PTDF_{ij,km}$ by:
\begin{equation}
\label{eq:flow}
    P_{ij} = \sum_{m} PTDF_{ij,km} P_m.
\end{equation}
The maximum power flow variations, in both directions, can then be evaluated as:
\begin{equation}
\label{eq:flow_max_pos}
    \Delta P_{ij}^{\text{max},+} = P_{ij}^{\text{max}} - P_{ij}
\end{equation}
\begin{equation}
\label{eq:flow_max_neg}
    \Delta P_{ij}^{\text{max},-} = - P_{ij}^{\text{max}} - P_{ij},
\end{equation}
where $\Delta P_{ij}^{\text{max},+}$ and $\Delta P_{ij}^{\text{max},-}$ are the maximum power flow variations respectively from $i$ to $j$ and from $j$ to $i$, and $P_{ij}^{\text{max}}$ is the line capacity.
Finally, we use that the change in the power flow of line $ij$ associated with a power injection at bus $m$ and equivalent withdrawal at $n$ can be obtained as:
\begin{equation}
    \Delta P_{ij} = (PTDF_{ij,kn}-PTDF_{ij,km}) \Delta P_{mn}.
    \label{eq:PlinePTDF}
\end{equation}
Algorithm \ref{algo:ptdf} describes how to evaluate the maximum quantity that can be traded for an injection in bus $m$ and retrieve in bus $n$.

\begin{algorithm}[t]
\small
\KwData{ \textit{request\_bus}, \textit{offer\_bus}, \textit{Quantity}}
\uIf{up regulation}{
$m$ = \textit{offer\_bus}\;
$n$ = \textit{request\_bus}\;
}
\uElseIf{down regulation}{
$m$ = \textit{request\_bus}\;
$n$ = \textit{offer\_bus}\;
}
\For{all the lines $ij$ in the distribution system}{
 Calculate $P_{ij}$ with (\ref{eq:flow}), $\Delta P_{ij}^{\text{max},+}$ with (\ref{eq:flow_max_pos}), $\Delta P_{ij}^{\text{max},-}$ with (\ref{eq:flow_max_neg})\;
 Calculate \textit{Quantity\_max} that can be injected in bus $m$ and retrieved in bus $n$ applying \cref{eq:PlinePTDF}, taking into account the direction of the flow\;
 Update \textit{Quantity} to be lower than or equal to  \textit{Quantity\_max};
 }
 \Return{Quantity}
 \caption{\small Calculation of the maximum quantity that can be exchanged between a request bus and an offer bus}
 \label{algo:ptdf}
\end{algorithm}
\normalsize

Using these equations, we can prove that, in the described setup, if the activation of the maximum capacity designated by the bid satisfies the network constraints, any partial activation of the bid will satisfy the network constraints as well.

Looking at \eqref{eq:PlinePTDF}, the term $PTDF_{ij,kn}-PTDF_{ij,km}=\alpha_{ij}$ is a constant, dependent only on the network topology and the line reactances. As a result, \eqref{eq:PlinePTDF} can become: 
\begin{equation}
\Delta \mathbf{P}_{\textnormal{line}}= \bm{\alpha} \Delta P_{mn},
\label{eq:proof1}
\end{equation} 
where $\bm{\alpha}=[\alpha_{ij}]$ is a vector of size $L \times 1$ with $L$ being the number of lines. Following the assumption in this paper that each bid is strictly either for up-regulation or down-regulation (but not both up and down in the same bid), a partial activation of the bid will be between 0 and $\Delta P_{mn}^\text{max}$ when the maximum capacity of the offered bid is activated. It follows that $\Delta P_{ij}^\text{max}= \alpha_{ij} \Delta P_{mn}^\text{max}$, and as a result the maximum change in each of the line flows occurs when the full bid, i.e., up to its maximum capacity, is activated. Assuming an up-regulation activation, as long as $|P_{ij}+\Delta P_{ij}^\text{max}| \leq P_{ij}^\text{max}$, for \emph{all} lines, it is straightforward to see that for any partial activation of the bid  $\Delta P_{mn} \leq \Delta P_{mn}^\text{max}$, it will hold $|P_{ij}+\Delta P_{ij}| \leq |P_{ij}+\Delta P_{ij}^\text{max}| \leq P_{ij}^\text{max}$. We can perform a similar derivation for any down-regulation bid.

\vspace{-0.1cm}
\subsection{Simulation and Results}
In this study, we show the organization of a market where several requests are submitted to the local flexibility market, including the location where the flexibility will be received. When there is a match with an offer in terms of price, the resulting potential power flows are evaluated to make sure that the activation would not lead to any congestion.
The bids used for this study are given in Tables \ref{tab:requests} and \ref{tab:results}. Without loss of generality, we assume that all requests are submitted as a batch, and offers are submitted one by one later. 
The corresponding market clearing, performed each time a bid is added, is described in Algorithm \ref{algo:clearing}.

\begin{algorithm}[t]
\small
\KwData{\textit{Power\_Dispatch, offer, All\_Requests}}
Compare \textit{offer} to all requests in the same direction: \\
\For{\textit{request} in \textit{All\_Requests}}{
Check that the prices match: \textit{offer\_price} $\leq$ \textit{request\_price}\;
Initialize \textit{Quantity} = min(\textit{offer\_quantity}, \textit{request\_quantity})\;
\For{c in all combinations of previously accepted requests}{
Modify \textit{Power\_Dispatch} to account for \textit{c} being activated\;
Calculate \textit{Quantity\_max} that can be exchanged between  \textit{offer\_bus} and \textit{request\_bus} for \textit{c}, applying Algorithm \ref{algo:ptdf}\;
Update \textit{Quantity} to be lower than or equal to \textit{Quantity\_max}\;
}
\If{Quantity $>$ 0}{
\textit{offer} and \textit{request} match for \textit{Quantity}\;
\If{request\_type is Unconditional}{
Update \textit{Power\_Dispatch} accordingly\;
}
Update and order \textit{Order\_Book}\;
}
}
\If{there was a match with an unconditional request}{
Try matching offers in \textit{Order\_Book} and repeat until no new match with an unconditional request is found.
}
 \Return{Order\_Book, Power\_Dispatch}
 \caption{\small Market clearing}
 \label{algo:clearing}
\end{algorithm}
\normalsize

 The results of the matching algorithm are shown in Table \ref{tab:results}. We can see that \textit{offer2} is partially matched with \textit{req3}, due to line congestion and the rest of the quantity in \textit{offer2} is added to the order book. Later, in bidding round number 5, the new \textit{offer5} can match with the rest of \textit{req3} without the risk of creating any congestion. On the other hand, \textit{offer3} cannot be procured because its activation could create a congestion.

\begin{table}[t]
\centering
\caption{Requests submitted to the market, by order of submission}
\label{tab:requests}
\begin{tabular}{lllccc}
   &  &            &  & Quantity & Price \\ 
ID   & Direction & Type           & Bus &  (kW) &  (€/kW) \\ \hline
req1 & Up        & Unconditional & 13  & 30     & 0.042    \\
req2 & Down      & Conditional   & 4   & 10     & 0.044    \\
req3 & Down      & Conditional   & 10  & 20     & 0.041    \\
req4 & Up        & Conditional   & 15  & 20     & 0.041    \\
req5 & Down      & Unconditional & 5   & 10     & 0.040    \\
req6 & Up        & Conditional   & 10  & 30     & 0.037    \\ \hline
\end{tabular}
\end{table}

\begin{table*}[t]
\centering
\caption{Offers submitted to the market, by order of submission, and matches with the requests given in Table \ref{tab:requests}}
\label{tab:results}
\begin{tabular}{llccc|llllll}
\multicolumn{5}{c}{Offers}     & \multicolumn{6}{|c}{Bidding rounds (quantities in kW)}                      \\ \hline
ID &
  Type &
  Bus &
  Quantity (kW) &
  Price (€/kW) &
  \multicolumn{1}{c}{1} &
  \multicolumn{1}{c}{2} &
  \multicolumn{1}{c}{3} &
  \multicolumn{1}{c}{4} &
  \multicolumn{1}{c}{5} &
  \multicolumn{1}{c}{6} \\ \hline
offer1 & Up   & 14 & 30 & 0.035 & req1: 30 &  &             &            &  &            \\
offer2 &
  Down &
  13 &
  40 &
  0.040 &
   &
  \begin{tabular}[c]{@{}l@{}}req2: 10\\ req3: 10\\ congestion\end{tabular} &
   &
   &
   &
   \\
offer3 & Down & 12 & 30 & 0.039 &            &  & congestion &            &  &            \\
offer4 & Up   & 15 & 20 & 0.032 &            &  &             & req4: 20 &  &            \\
offer5 &
  Down &
  8 &
  40 &
  0.033 &
   &
   &
   &
   &
  \begin{tabular}[c]{@{}l@{}}req3: 10\\ req5 10\end{tabular} &
   \\
offer6 & Up   & 7  & 40 & 0.031 &            &  &             &            &  & req6: 30 \\ \hline
\end{tabular}
\vspace{-0.2cm}
\end{table*}

 \vspace{-0.1cm}
\section{Discussion}\label{sec:Disc}
This section discusses some of the market design characteristics introduced in this paper, their limitations, and possible alternatives, which we hope will inspire researchers for future work in the wider area of local flexibility markets that include network constraints.
\paragraph{Power Flow Algorithms} \label{sec:DiscPF}
This paper has used a DC power flow algorithm for the network check as a first step towards the inclusion of network constraints in a continuous market-clearing algorithm. The motivation behind this choice is that (i) it is simpler, and, thus, it is more transparent for the market players while it also allows us to obtain valuable insights about the market design choices we had to make, (ii) it is a linear algorithm, and as such it is faster to solve as it does not require iterations, and (iii) it allows the use of PTDFs, which enable us to extract a useful proof when it comes to assessing the impact of partial versus full activation of matched bids. Although we consider the DC power flow algorithm a valuable first step to include active power flows and consider line congestion, it also introduces two main limitations. First, distribution lines are not characterized by a much higher reactance compared to the ohmic resistance, and therefore the DC power flow approximation might misestimate the actual line flows. Second, voltage issues are more common in the distribution grid, and reactive power flows need to be considered. Similarly, line losses can be non-negligible. Therefore, future versions of such markets shall investigate the application of power flow algorithms such as \textit{LinDistFlow} or \textit{DistFlow} \cite{BaranWu1989_LinDistFlow}
, or even AC Power Flow \cite{notes}. 

\paragraph{Active and Reactive Power Markets} Considering that distribution grids often face voltage issues that could be resolved through appropriate local control of reactive power, local flexibility markets can offer an ideal platform for trading reactive power. Therefore, we think that an extension of the proposed market setup to include a reactive power market would be valuable to be explored. 

\paragraph{Estimation of Baseline Dispatch and Location of FlexRequests} In this paper, we suggest that the FMO has knowledge of the baseline dispatch either because the energy markets have cleared before the flexibility market, or the DSO has offered its best estimate, or the FMO has collected data and determined the most probable scenarios. Similarly, in this market setting, we suggest that \textit{FlexRequests} shall determine the location or, otherwise, the FMO should estimate a set of probable locations for each \textit{FlexRequest} and assess their impact on network violations. In both cases, scenarios need to be assumed. Considering that it is impossible to account for all possible scenarios, an extension of the proposed market clearing algorithm is to formulate it as a probabilistic market clearing, and allocate a small amount of reserves to counter any instance not captured by the scenarios.

\paragraph{Multi-Period Market Clearing} \textit{FlexOffers} by distributed resources often have a rebound effect. For instance, thermal loads and batteries have to replenish the energy they offered at a later point in time. Therefore, a local flexibility market should consider block offers and bids that span multiple time periods. This will be object of our future work.

\paragraph{Integration with Existing Markets} The market design we propose in this paper can consider both energy and reserves (unconditional and conditional offers) and is suitable for any time resolution (e.g., month-ahead, week-ahead, day-ahead, intra-day). Future work shall focus on ways that could optimally integrate such a local flexibility market to the existing energy and reserve markets both at the wholesale level, and in the future distribution-level markets.

\paragraph{Market Power} One of the key criteria for any market design is the (in)ability of the market players to exert market power. In a local flexibility market, the bidding strategies could depend on whether the market players receive a price only for the power they offer, or also for the offered energy during activation. Detailed analyses about the potential to exert market power need to be carried out to investigate potential issues and compare them with alternative designs.

\paragraph{Activation of \textit{FlexOffers} in Real-Time} In this paper, we suggest that although BRPs and the DSO compete for flexibility reserves, during real-time it can help avoid any possible network violations if the DSO is able to activate some of the \textit{FlexOffers} procured, on top of the ones activated. This ensures network feasibility. Alternative directions to address this challenge could also be sought. 

\paragraph{Matching Up- and Down-Regulation Bids} This paper suggests to separate the up-regulation from the down-regulation bids and treat them individually and separately. This allows for higher flexibility, as certain resources may be able to (or prefer to) offer only up-regulation or down-regulation (e.g., solar PVs). At the same time, we require that each \emph{FlexRequest} for up-regulation is matched with a \emph{FlexOffer} for up-regulation, and similarly for down-regulation. However, network constraints add a new dimension of complexity. Cases can exist that a given \emph{FlexRequest} could be better served (i.e., cheaper) by a mix of \emph{FlexOffers} for up- and down-regulation at different locations of the grid. We plan to look into such extension of our proposed market design in our future work. 

\paragraph{Pricing Mechanisms} In this paper, we have assumed that each market participant pays or gets paid the price they have bid. For bids to match, the \textit{FlexRequest} price shall either be higher or equal to the price of the \textit{FlexOffer}. Here, we follow the common approach for pricing: in every matched pair, the price for both \textit{FlexRequest} and \textit{FlexOffer} is equal to the price of the first incoming bid. However, alternative pricing mechanisms can also exist. The price difference could be allocated to the FMO, and then be used for network investments by the DSO or other purposes determined by the regulator. Alternatively, one could set the price equal to the lowest price of the two (i.e., the \textit{FlexOffer}). Further assessment related to the implications on social welfare and bidding strategies is required to identify the most appropriate pricing mechanisms for different local flexibility markets.


 \vspace{-0.1cm}
\section{Conclusion}\label{sec:Conc}
To the best of our knowledge, this paper proposes for the first time a design of a \emph{continuous} local flexibility market that explicitly considers network constraints. We discuss the general architecture of such a market, the structure of the \textit{FlexRequests}, and elaborate on a number of design options for the inclusion of network constraints in the market clearing. In the early stages of  local flexibility markets, where insufficient liquidity may hinder market development, continuous markets are expected to be the most suitable option. At the same time, in increasingly loaded distribution systems, including the network constraints in the market clearing ensures that every matched pair of bids will not violate operational limits, and would not require additional actions from the distribution system operators that result in additional costs. This paper focuses on active power markets and it has integrated linearized power flow equations (DC power flow) to ensure no line limit violations. Aiming to establish a starting point for future research on specific design parameters of local flexibility markets, in the last part of this paper, we discuss a series of questions and research challenges that require further exploration and assessment. In our future work, we intend to include a higher level of detail of the power flow equations and establish a common framework for active and reactive power local flexibility markets.

 \vspace{-0.1cm}
\bibliographystyle{IEEEtran}

\bibliography{Bib}

\end{document}